\documentclass[11pt,a4paper]{article}
\usepackage[margin=30mm,top=27mm,bottom=28mm,headheight=14pt]{geometry}
\usepackage{amsmath,amssymb,amsthm}
\usepackage{fancyhdr}
\usepackage[hidelinks,bookmarksnumbered]{hyperref}

\hypersetup{
  pdftitle={Unbounded derived length for finite p-groups with a splitting automorphism},
  pdfauthor={Jiang Yang},
  pdfsubject={Finite groups and splitting automorphisms of prime-power order},
  pdfkeywords={finite p-group, splitting automorphism, derived length, wreath product}
}
\newtheorem{theorem}{Theorem}[section]
\newtheorem{lemma}[theorem]{Lemma}
\DeclareMathOperator{\dl}{dl}
\DeclareMathOperator{\Sym}{Sym}
\newcommand{\ev}{\operatorname{ev}}
\newcommand{\Z}{\mathbb Z}
\newcommand{\dd}{d}
\newcommand{\doi}[2]{\href{https://doi.org/#1}{#2}}
\title{\bfseries Unbounded derived length for finite \(p\)-groups\\
       with a splitting automorphism}
\author{Jiang Yang\\[4pt]
  \small School of Mathematical Sciences, Guangxi Minzu University\\
  \small Nanning 530006, China\\
  \small\href{mailto:yangjiangdy@126.com}{\texttt{yangjiangdy@126.com}}}
\date{}

\begin{document}
\maketitle

\begin{abstract}
For every integer \(\ell\ge 2\) there is a two-generator finite \(2\)-group
of derived length \(\ell\) admitting a splitting automorphism of exact
order eight. This gives a negative answer to Problem~15.96(a) of the
\emph{Kourovka Notebook}. More generally, let \(T\) be a nontrivial finite
\(p\)-group of exponent \(m\). We construct a finite \(p\)-group \(B\)
requiring exactly \(p\) generators, with a splitting automorphism of exact
order \(pm\), such that \(B'\) is a subdirect product of copies of \(T\).
In particular, \(\dl(B)=1+\dl(T)\). For each fixed prime power
\(m=p^e\ge4\), this yields unbounded derived length at the fixed
automorphism order \(pm\).
\end{abstract}


\section{Introduction}\label{sec:intro}

An automorphism \(\alpha\) of a group \(G\) is \emph{splitting of order
\(q\)} in the sense of the \emph{Kourovka Notebook} if \(\alpha^q=1\) and
\begin{equation}\label{eq:norm}
  g\,\alpha(g)\,\alpha^2(g)\cdots\alpha^{q-1}(g)=1
  \qquad(g\in G).
\end{equation}
Khukhro asked whether the derived length of a \(d\)-generator nilpotent
\(p\)-group admitting such an automorphism for \(q=p^n\) is bounded in
terms of \(d,p,n\)~\cite[Problem~15.96(a), p.~99]{notebook}. For \(n=1\),
Khukhro proved that even the nilpotency class is bounded in terms of
\(d\) and \(p\)~\cite{khukhro}. Ersoy proved that a finite group
admitting a splitting automorphism of odd order is soluble~\cite{ersoy}.

We write $d(G)$ for the least number of generators of a group $G$
and $G'=[G,G]$ for its derived subgroup.
For a finite group $G$, its exponent $\exp(G)$ is the least positive
integer $m$ such that $g^m=1$ for every $g\in G$.
For a soluble group $G$, we write $\operatorname{dl}(G)$ for its
derived length, with $\operatorname{dl}(1)=0$.
For background on group theory, see \cite{Robinson1996}.

\begin{theorem}\label{thm:main}
Let \(p\) be a prime and let \(m\ge 4\) be a power of \(p\). For every
\(\ell\ge 2\) there is a finite \(p\)-group \(B\) such that
\[
    \dd(B)=p,\qquad \dl(B)=\ell,\qquad \exp(B')=m,
\]
and \(B\) admits a splitting automorphism of exact order \(pm\).
\end{theorem}

For \(p=2\) and \(m=4\), the parameters in Khukhro's question are
\((d,p,n)=(2,2,3)\) for every \(\ell\). The proposed bound in
part~(a) therefore fails already for two generators and automorphism
order eight. Part~(b) of the same problem concerns order four.

\pagebreak[3]
Razmyslov's theorem gives finite \(p\)-groups of fixed prime-power
exponent \(m\ge4\) and arbitrarily large derived length
\cite[Chapter~4]{vaughanlee}. Razmyslov settled the
exponent-four case~\cite{razmyslov}, and Bachmuth, Mochizuki and Walkup
had earlier treated exponent five~\cite{bmw}. De Cornulier and
Mann used this theorem to construct two-generator residually \(p\)-groups
whose derived subgroups have bounded exponent and whose finite
\(p\)-quotients have unbounded derived length~\cite{cornuliermann}.
In our construction the \(p\)-th powers of all elements of a single
coset lie in a base group of exponent \(m\). The values of a function
on two permutation orbits show that the \(m\)-th power of the induced
automorphism moves a translation.

\section{A finite wreath product construction}\label{sec:construction}

Let \(A\) be a finite abelian group, written additively. Write
\(\tau_a(x)=a+x\) for its regular translations and \(\tau(A)\) for
the group they form. For a group \(T\), write
\(T^A\) for the group of all functions from \(A\) to \(T\), with
pointwise multiplication. We use the action
\[
  (\pi\cdot f)(x)=f(\pi^{-1}x)
  \qquad(\pi\in\Sym(A),\ f\in T^A)
\]
in \(T^A\rtimes\Sym(A)\). We identify a permutation with its natural
lift to this semidirect product. All conjugations below are left
conjugations.

\begin{lemma}\label{lem:coordinate}
Suppose that \(f:A\to T\) has a value equal to \(1\) and that its
values generate \(T\). Let \(B=\langle f,\tau(A)\rangle\).
Then \(B'\le T^A\), and evaluation at each \(x\in A\) maps \(B'\)
onto \(T\).
\end{lemma}

\begin{proof}
The permutation image of \(B\) is abelian, so \(B'\le T^A\).
Write \(f_a=\tau_a f\tau_a^{-1}\). For \(a,b\in A\), the images of
\(f_a\) and \(f_b\) in \(B/B'\) are both the image of \(f\). Hence
\(f_af_b^{-1}\in B'\), and
\[
    \ev_x(f_af_b^{-1})=f(x-a)f(x-b)^{-1}.
\]
Choose \(b\) so that \(f(x-b)=1\). As \(a\) varies, the displayed
values include every value of \(f\). They therefore generate \(T\).
\end{proof}

\begin{theorem}\label{thm:transfer}
Let \(T\) be a nontrivial finite \(p\)-group of exponent \(m\).
There is a finite \(p\)-group \(B\) with \(\dd(B)=p\) and a splitting
automorphism \(\alpha\) of exact order \(pm\). The image of
\(\alpha\) in \(\operatorname{Out}(B)\) has order \(p\), and \(B'\)
is a subdirect product of finitely many copies of \(T\). In particular,
\[
       \dl(B)=1+\dl(T),\qquad \exp(B')=m .
\]
\end{theorem}

\begin{proof}
Choose \(v\in T\) of order \(m\) and generators \(g_1,\ldots,g_k\)
of \(T\). Such a \(v\) exists since all element orders are powers of
\(p\). Choose \(a\ge2\) with \(p^{a(p-1)}>2p^2+pk\), let
\(R=\Z/p^a\Z\), and put
\[
 A=\{(x_0,\ldots,x_{p-1})\in R^p:\ \textstyle\sum_{i=0}^{p-1}x_i=0\}.
\]
This is a free \(R\)-module of rank \(p-1\). Let \(\rho\) cyclically
permute the coordinates, and put \(u=e_0-e_1\), where
\(e_0,\ldots,e_{p-1}\) are the coordinate vectors of \(R^p\).
The restriction of \(\rho\) to \(A\) has order \(p\), since
\(\rho(u)\ne u\).
For every \(b\in A\), each coordinate of
\(b+\rho b+\cdots+\rho^{p-1}b\) is the sum of the coordinates of \(b\).
Consequently
\begin{equation}\label{eq:affine}
      (\tau_b\rho)^p=\tau_{b+\rho b+\cdots+\rho^{p-1}b}=1
      \qquad(b\in A).
\end{equation}
In particular, \(s=\tau_u\rho\) has order \(p\), since \(s(0)=u\ne0\).
Its orbit \(O_0\) through \(0\) has length \(p\).

\pagebreak[2]
Since \(|A|>2p^2+pk\), we can take a second \(\langle s\rangle\)-orbit
\(O_1\). We construct a \(\rho\)-invariant
function \(f:A\to T\). Set \(f(0)=v\). Give the value \(1\) to
every \(\langle\rho\rangle\)-orbit meeting
\((O_0\setminus\{0\})\cup O_1\). These assignments agree because
\(\rho\) fixes \(0\). The \(\langle\rho\rangle\)-orbits meeting
\(O_0\cup O_1\) use at most \(2p^2\) points, and each remaining
orbit has at most \(p\) points. There are therefore \(k\) distinct
unused orbits on which to place \(g_1,\ldots,g_k\), respectively.
Set \(f=1\) on all other orbits.

Inside \(T^A\rtimes\langle\tau(A),\rho\rangle\), put
\[
     B=\langle\tau(A),f\rangle,\qquad t=f s .
\]
The ambient group is a finite \(p\)-group. The element \(\rho\)
normalizes the translations and fixes \(f\), so it normalizes \(B\).
Since \(f,\tau_u\in B\), the element \(t=f\tau_u\rho\) normalizes
\(B\) as well. Write \(\alpha(g)=tgt^{-1}\) on \(B\).

For every \(g\in B\), the permutation part of \(gt\) is
\(\tau_b\rho\) for some \(b\in A\). Equation~\eqref{eq:affine} gives
\((gt)^p\in T^A\). The base group \(T^A\) has exponent \(m\), so
\((gt)^{pm}=1\) for every \(g\in B\), including \(g=1\).
Putting \(q=pm\), we have \(t^q=1\), so \(\alpha^q=1\).
The ordered product in~\eqref{eq:norm} is trivial by the identity
\begin{equation}\label{eq:expand}
    (gt)^q=\bigl(g\,\alpha(g)\cdots\alpha^{q-1}(g)\bigr)t^q .
\end{equation}

It remains to determine the order of \(\alpha\). The base function
\(t^p\) has the coordinate formula
\[
    t^p(x)=f(x)f(s^{-1}x)\cdots f(s^{-(p-1)}x).
\]
On \(O_0\) the factors include \(v\) exactly once and are otherwise
equal to \(1\). Every factor on \(O_1\) equals \(1\). Since
\(t^m=(t^p)^{m/p}\), we have
\begin{equation}\label{eq:marker}
    (t^m)(x)=v^{m/p}\quad(x\in O_0),\qquad
    (t^m)(x)=1\quad(x\in O_1).
\end{equation}
The element \(v^{m/p}\) has order \(p\), so this base function
is nonconstant. A base function commuting with every
translation is constant, by transitivity of the translations. Thus
\(t^m\) fails to centralize \(B\), and \(\alpha^m\ne1\). Together with
\(\alpha^{pm}=1\) and the fact that \(pm\) is a power of \(p\), this
gives \(|\alpha|=pm\). Moreover, \(t\in B\rho\) and \(\rho^p=1\),
so \(t^p\in B\) and \(\alpha^p\) is inner. On the abelian quotient
\(B/(B\cap T^A)\cong A\), the induced action is \(\rho\), of order
\(p\). Thus the outer automorphism represented by \(\alpha\) has order \(p\).

The translations associated to an \(R\)-basis of \(A\), together with
\(f\), generate \(B\). Omitting a basis translation gives a proper
permutation image, and omitting \(f\) gives just the translation group.
This generating set is therefore irredundant. The Burnside basis
theorem now gives \(\dd(B)=p\). Lemma~\ref{lem:coordinate} applies
to \(f\), and hence \(B'\) is a subdirect product of copies of \(T\).
The embedding in \(T^A\) gives \(\dl(B')\le\dl(T)\) and
\(\exp(B')\mid m\). Projection onto a factor gives the reverse
inequality and divisibility. As \(T\ne1\), it follows that
\(\dl(B)=1+\dl(T)\).
\end{proof}

\pagebreak[3]
\section{The counterexamples}\label{sec:examples}

\begin{proof}[Proof of Theorem~\ref{thm:main}]
Put \(r=\ell-1\). By Razmyslov's theorem there is a locally finite
nonsoluble group \(L\) of exponent \(m\)~\cite[Chapter~4]{vaughanlee}.
The derived word \(\delta_{r-1}\), with \(\delta_0(x)=x\), has a
nontrivial value in \(L\). Its entries generate a finite \(p\)-group
\(H\) of exponent dividing \(m\) with \(H^{(r-1)}\ne1\), where
\(H^{(i)}\) is the \(i\)-th derived subgroup and \(H^{(0)}=H\).
See also the finite formulation in~\cite[Theorem~4]{cornuliermann}.
The group \(H/H^{(r)}\) has derived length exactly
\(r\), since the derived series of a finite soluble group is strictly
descending until it reaches \(1\). Thus
\[
        T=(H/H^{(r)})\times C_m
\]
has derived length \(r\) and exponent \(m\). Theorem~\ref{thm:transfer}
applied to this \(T\) gives all three assertions about \(B\) and the
required automorphism.
\end{proof}

For order eight one can see the order calculation on four points.
Here \(A\cong\Z/2^a\Z\), \(\rho(x)=-x\) and \(s(x)=1-x\).
The pairs \(\{0,1\}\) and \(\{2,-1\}\) are \(s\)-orbits.
Take \(f(0)=v\) of order four and \(f(\pm1)=f(\pm2)=1\), putting the
generators of \(T\) on other pairs \(\{x,-x\}\).
Then \(B=\langle \tau_1,f\rangle\), and \(t^4\) takes the values
\(v^2\) and \(1\) on the two indicated pairs. In particular \(t^4\)
does not commute with \(\tau_1\), although every element of \(Bt\)
has eighth power \(1\).

\paragraph{Methods.}
The author conceived the affine-action approach used in Theorem~2.2.
Open AI's ChatGPT (Codex) assisted the author with the calculations
needed to carry out the construction, including those involved in
selecting the two orbits and checking the coordinate projections,
the ordered product, and the order of the induced automorphism.

\subsection*{Declaration of generative AI and AI-assisted technologies\\
in the manuscript preparation process}
The author used Open AI's ChatGPT (Codex) to help draft and revise
the English text. The author takes
responsibility for the mathematical content and the final manuscript.


\begin{thebibliography}{9}\small
\bibitem{bmw}
S.~Bachmuth, H.~Y.~Mochizuki and D.~Walkup,
\doi{10.1090/S0002-9904-1970-12469-7}{A nonsolvable group of exponent \(5\)},
\emph{Bull. Amer. Math. Soc.} \textbf{76} (1970), 638--640.

\bibitem{cornuliermann}
Y.~de Cornulier and A.~Mann,
\doi{10.1007/978-3-7643-8412-8_4}{Some residually finite groups satisfying laws},
in \emph{Geometric Group Theory}, Trends in Mathematics, Birkh\"auser,
Basel, 2007, 45--50.

\bibitem{ersoy}
K.~Ersoy,
\doi{10.1007/s00013-016-0874-6}{Finite groups with a splitting automorphism of odd order},
\emph{Arch. Math. (Basel)} \textbf{106} (2016), 401--407.

\bibitem{khukhro}
E.~I.~Khukhro,
\doi{10.1070/SM1987v058n01ABEH003095}{On locally nilpotent groups admitting a splitting automorphism of prime order},
\emph{Math. USSR-Sb.} \textbf{58} (1987), 119--126.

\bibitem{notebook}
E.~I.~Khukhro and V.~D.~Mazurov (eds.),
\href{https://arxiv.org/abs/1401.0300v50}{\emph{Unsolved Problems in Group Theory. The Kourovka Notebook}},
21st ed., Sobolev Institute of Mathematics, Novosibirsk, 2026,
November revision, arXiv:1401.0300v50.

\bibitem{razmyslov}
Yu.~P.~Razmyslov,
\doi{10.1070/IM1979v013n01ABEH002015}{On a problem of Hall and Higman},
\emph{Math. USSR-Izv.} \textbf{13} (1979), 133--146.

\bibitem{Robinson1996}
D. J. S. Robinson,
\emph{A Course in the Theory of Groups},
2nd ed., Graduate Texts in Mathematics, vol.~80,
Springer-Verlag, New York, 1996.

\bibitem{vaughanlee}
M.~Vaughan-Lee, \emph{The Restricted Burnside Problem}, 2nd ed.,
London Mathematical Society Monographs, New Series, vol.~8, Clarendon
Press, Oxford, 1993, \doi{10.1093/oso/9780198537861.003.0004}{Chapter~4, 67--86}.
\end{thebibliography}
\end{document}